\documentclass[12pt]{amsart}
\usepackage[utf8]{inputenc}
\usepackage{graphicx}
\usepackage{epstopdf}
\usepackage{amsmath}
\usepackage{xcolor}
\usepackage{amssymb}
\usepackage{tikz}
\usetikzlibrary{arrows}

\newtheorem{question}{Question}
\newtheorem*{ack}{Acknowledgements}

\makeatletter
\def\blfootnote{\gdef\@thefnmark{}\@footnotetext}
\makeatother

\def\house#1{\setbox1=\hbox{$\,#1\,$}%
\dimen1=\ht1 \advance\dimen1 by 2pt \dimen2=\dp1 \advance\dimen2 by 2pt
\setbox1=\hbox{\vrule height\dimen1 depth\dimen2\box1\vrule}%
\setbox1=\vbox{\hrule\box1}%
\advance\dimen1 by .4pt \ht1=\dimen1
\advance\dimen2 by .4pt \dp1=\dimen2 \box1\relax}

\begin{document}
\title{Representation of similar triangles in $\mathbb{R}^{3}$}
\author{Jaitra Chattopadhyay and Siddhartha Sankar Chattopadhyay}
   \address[Jaitra Chattopadhyay]{Department of Mathematics, Siksha Bhavana, Visva-Bharati\\ Santiniketan - 731235, West Bengal, India}
\email{jaitra.chattopadhyay@visva-bharati.ac.in; chat.jaitra@gmail.com}
\address[Siddhartha Sankar Chattopadhyay]{Mathematics Teacher (Retd.)\\
Bidhannagar Govt. High School, BD-303\\
Salt Lake, Kolkata - 700064\\
West Bengal, India}
\email{1959ssc@gmail.com}

\begin{abstract}
In 1994, Dhar and Sinha provided a graphic representation in $\mathbb{R}^{2}$ of the classification of triangles based on similarity. In this article, we briefly recall their treatment and offer a generalization by providing a graphic representation in $\mathbb{R}^{3}$. Towards the end of this article, we make an attempt to offer a quantitative analysis of some particular class of triangles.
\end{abstract}
\maketitle

\section{Introduction}

Let $\Delta$ be a triangle with side lengths $a,b$ and $c$ units. If $\Delta^{\prime}$ is another triangle with side lengths $a^{\prime},b^{\prime}$ and $c^{\prime}$ units such that $$\frac{a}{a^{\prime}} = \frac{b}{b^{\prime}} = \frac{c}{c^{\prime}},$$ then we say that $\Delta$ is similar to $\Delta^{\prime}$ and denote it by $\Delta \sim \Delta^{\prime}$. It is quite straightforward to check that the relation $`\sim'$ on the set of all triangles is an equivalence relation and therefore, this induces a partition of this set into certain equivalence classes. Thus each triangle in the euclidean plane $\mathbb{R}^{2}$ belongs to exactly one equivalence class. We define $[\Delta]$ to be $$[\Delta] = \{\Delta^{\prime} : \Delta^{\prime} \mbox{ is a triangle and } \Delta \sim \Delta^{\prime}\}$$ and denote the set of all equivalence classes of triangles by $\mathcal{T}$.

\smallskip

In \cite{ds}, Dhar and Sinha beautifully classified triangles based on similarity and provided a pictorial representation. Their article appeared in Indian Journal of Mathematics Teaching and due to its unavailability on the web, here we briefly recall their treatment before going into our method of representation.  

\section{Representation in $\mathbb{R}^{2}$ done by Dhar and Sinha}\label{ds-section}

Let the sides of a triangle $\Delta$ in an equivalence class based on similarity be $a,b$ and $c$. Without loss of any generality, let us assume that $a \leq b \leq c$ and let $x = \frac{a}{c}$ and $y = \frac{b}{c}$. Then we have $0 < x \leq 1$, $0 < y \leq 1$, $x \leq y$ and by virtue of the triangle inequality, we have $x + y > 1$. Again, for any $\Delta^{\prime} \in [\Delta]$ with corresponding sides $a^{\prime},b^{\prime}$ and $c^{\prime}$, we have $a^{\prime} \leq b^{\prime} \leq c^{\prime}$, $\frac{a^{\prime}}{c^{\prime}} = \frac{a}{c} = x$ and $\frac{b^{\prime}}{c^{\prime}} = \frac{b}{c} = y$. In other words, the values of $x$ and $y$ are of constant values for all triangles in $[\Delta]$. 

\smallskip

This allows us to define a map from $\mathcal{T}$ to $\mathbb{R}^{2}$ by sending an equivalence class of triangles with sides $a,b,c$ with $a \leq b \leq c$ to the point $(\frac{a}{c},\frac{b}{c}) \in \mathbb{R}^{2}$. We see that the image of $\mathcal{T}$ under this map is the region enclosed by the red lines depicted in Figure 1. We observe that the points $(x,y)$ on the line segment $BE$ satisfy $x + y = 1$ and therefore do not lie on the image of $\mathcal{T}$. 

\smallskip

Next, Dhar and Sinha considered the equivalence classes of equilateral, isosceles, right-angled triangles and investigated which subsets of the aforementioned region correspond to each class. For instance, if $\Delta$ is equilateral, then we have $a = b = c$ and thus $x = y = 1$.

\begin{center}

\usetikzlibrary{calc}
\begin{tikzpicture}[x=0.5cm,y=0.5cm,z=0.3cm,>=stealth]
\draw[->] (xyz cs:x=0) -- (xyz cs:x=13.5) node[above] {$x$};
\draw[->] (xyz cs:y=0) -- (xyz cs:y=13.5) node[right] {$y$};

\node[fill,circle,inner sep=1.5pt,label={below:$A(1,0)$}] at (10,0) {};
\node[fill,circle,inner sep=1.5pt,label={left:$B(0,1)$}] at (0,10) {};
\node[fill,circle,inner sep=1.5pt,label={right:$C(1,1)$}] at (10,10) {};
\node[fill,circle,inner sep=1.5pt,label={left:$E(\frac{1}{2},\frac{1}{2})$}] at (5,5) {};
\node[fill,circle,inner sep=1.5pt,label={below:$O(0,0)$}] at (0,0) {};
\node[fill,circle,inner sep=1.5pt,label={left:$D(\frac{1}{\sqrt{2}},\frac{1}{\sqrt{2}})$}] at (7.1,7.1) {};

\draw [solid] (10,0) -- (10,10);
\draw [red] (0,10) -- (10,10);
\draw [solid] (0,0) -- (5,5);
\draw [solid] [red] (5,5) -- (10,10);
\draw [dashed] [red] (0,10) -- (5,5);
\draw [dashed] (5,5) -- (10,0);
\draw [solid] (0,0) -- (10,0) arc[start angle = 0, end angle = 90, radius = 5cm] -- (0,0);

\end{tikzpicture}

\end{center}
\begin{center}
Figure 1
\end{center}

Thus the point $C$ in Figure 1 corresponds to the equivalence class of all equilateral triangles. Similarly, if we are interested in right-angled triangles, then we must have $a^{2} + b^{2} = c^{2}$, thanks to Pythagoras' theorem. This implies that $x^{2} + y^{2} = 1$ and hence all the points on the circular arc $BD$, excluding $B$, correspond to the equivalence classes of right-angled triangles. The point of intersection of the unit circle $\{(x,y) : x^{2} + y^{2} = 1\}$ and the straight line $y = x$ in the first quadrant is $D(\frac{1}{\sqrt{2}},\frac{1}{\sqrt{2}})$ which corresponds to the equivalence class of isosceles right-angled triangles. Likewise, the region below the arc $BD$ and above the straight line $BE$ stands for the equivalence class of obtuse-angled triangles and the region above the arc $BD$ signifies the equivalence class of acute-angled triangles. Therefore, the arc $BD$ can be viewed as a demarcation between the classes of acute-angled and obtuse-angled triangles. More interestingly, we note that for each point $(x,y)$ on the line segment $BC$, we have $y = 1$. Therefore, the points on $BC$, except $B$, correspond to the equivalence classes of triangles with side lengths $a,b,c$ satisfying $b = c$. In other words, except $C$, they represent the equivalence classes of acute-angled isosceles triangles whose equal sides are bigger than the third. Similarly, the points on $CD$, except $C$ and $D$, represent the equivalence classes of acute-angled isosceles triangles whose equal sides are smaller than the third side. Finally, the points on $DE$ except $D$ and $E$, represent the equivalence classes of obtuse-angled isosceles triangles.

\smallskip

From the above discussion about the work of Dhar and Sinha, we obtain a graphic representation of triangles in a specified region in the euclidean plane $\mathbb{R}^{2}$. It is worthwhile to note that this representation of similar triangles uses information only about the sides of a given triangle and it does not take angles into consideration. Even when representing the region corresponding to right-angled triangles, it uses Pythagoras' theorem which is an equation involving only the sides. A natural question at this point might be the following.
\begin{question}\label{quest}
How does one represent different equivalence classes of triangles by taking the angles into consideration?
\end{question}

\section{Representation of triangles in $\mathbb{R}^{3}$}\label{amader}

The main theme of this article is to probe into Question \ref{quest} and offer another representation similar to that by Dhar and Sinha. To proceed to answer Question \ref{quest}, we denote the angles of a triangle by the letters $x,y$ and $z$ to keep the exposition notationally convenient. We know that if $x,y,z$ denote the angles of a triangle, then we have $x + y + z = \pi$. We remark here that we do not distinguish between different permutations of the triplet $(x,y,z)$. 

\smallskip

Now, if $\Delta$ is a triangle with angles $x,y,z$, then $0 < x,y,z < \pi$ and $x + y + z = \pi$. We can therefore, define a map $\Psi$ from $\mathcal{T}$ to $\mathbb{R}^{3}$ by sending $[\Delta]$ to $(x,y,z) \in \mathbb{R}^{3}$. Therefore, the image of $[\Delta]$ under $\Psi$ is the portion of the plane $\Sigma = \{(x,y,z) \in \mathbb{R}^{3} : x + y + z = \pi\}$ lying in the first octant. We refer the reader to Figure 2 where the straight lines joining $A,B$ and $C$ consecutively are marked in dotted lines to indicate that the points of $\mathbb{R}^{3}$ lying on these lines do not correspond to any triangle since that would mean one of $x,y$ or $z$ is $0$ which is not admissible for a triangle. 

\smallskip

Now, we consider different equivalence classes of triangles and see which subsets of the plane $\Sigma$ correspond to those classes.

\bigskip

\noindent
{\bf Equilateral triangles:-} In this case, we have $x = y = z = \frac{\pi}{3}$ and therefore, the equivalence class of all equilateral triangles is mapped to the singleton $\{(\frac{\pi}{3},\frac{\pi}{3},\frac{\pi}{3})\} \in \Sigma$. 

\begin{center}

\begin{tikzpicture}[x=0.5cm,y=0.5cm,z=0.3cm,>=stealth]
\draw[->] (xyz cs:x=0) -- (xyz cs:x=13.5) node[below] {$x$};
\draw[->] (xyz cs:y=0) -- (xyz cs:y=13.5) node[right] {$z$};
\draw[->] (xyz cs:z=0) -- (xyz cs:z=-13.5) node[above] {$y$};

\node[fill,circle,inner sep=1.5pt,label={left:$P(0,\frac{\pi}{2},\frac{\pi}{2})$}] at (0,5,-5) {};
\node[fill,circle,inner sep=1.5pt,label={right:$Q(\frac{\pi}{2},0,\frac{\pi}{2})$}] at (5,5,0) {};
\node[fill,circle,inner sep=1.5pt,label={below:$R(\frac{\pi}{2},\frac{\pi}{2},0)$}] at (5,0,-5) {};
\node[fill,circle,inner sep=1.5pt,label={left:$A(0,0,\pi)$}] at (0,10,0) {};
\node[fill,circle,inner sep=1.5pt,label={above right:$B(\pi,0,0)$}] at (10,0,0) {};
\node[fill,circle,inner sep=1.5pt,label={below right:$C(0,\pi,0)$}] at (0,0,-10) {};
\node[fill,circle,inner sep=1.5pt,label={above right:$O$}] at (0,0,0) {};
\node[fill,circle,inner sep=1.5pt,label={right:$S$}] at (5,5,-5) {};
\node[fill,circle,inner sep=1.5pt,label={above right:$D(0,\frac{\pi}{2},0)$}] at (0,0,-5) {};
\node[fill,circle,inner sep=1.5pt,label={above right:$E(0,0,\frac{\pi}{2})$}] at (0,5,0) {};
\node[fill,circle,inner sep=1.5pt,label={above right:$F(\frac{\pi}{2},0,0)$}] at (5,0,0) {};
\node[fill,circle,inner sep=1.5pt,label={above right:$R^{\prime}$}] at (2.5,5,-2.5) {};
\node[fill,circle,inner sep=1.5pt,label={right:$Q^{\prime}$}] at (2.5,2.5,-5) {};
\node[fill,circle,inner sep=1.5pt,label={above right:$P^{\prime}$}] at (5,2.5,-2.5) {};

\draw [dashed] [red] (10,0,0) -- (0,10,0);
\draw [dashed] [red] (10,0,0) -- (0,0,-10);
\draw [dashed] [red] (0,0,-10) -- (0,10,0);
\draw [solid] (5,0,0) -- (5,5,0);
\draw [solid] (0,0,-5) -- (0,5,-5);
\draw [solid] (0,0,-5) -- (5,0,-5);
\draw [solid] (5,0,-5) -- (5,0,0);
\draw [solid] (5,0,-5) -- (5,5,-5);
\draw [solid] (5,5,-5) -- (0,5,-5);
\draw [solid] (5,5,-5) -- (5,5,0);
\draw [solid] (0,5,-5) -- (0,5,0); 
\draw [solid] (0,5,0) -- (5,5,0); 
\draw [blue] (0,5,-5) -- (5,0,-5);
\draw [blue] (5,0,-5) -- (5,5,0);
\draw [blue] (0,5,-5) -- (5,5,0);

%
%
\end{tikzpicture}

\end{center}
\begin{center}
Figure 2
\end{center}

\bigskip

\noindent
{\bf Right-angled triangles:-} In this case, one of $x$, $y$, $z$ equals $\frac{\pi}{2}$. Assume first that $z = \frac{\pi}{2}$. Then we have $x + y = \pi - \frac{\pi}{2} = \frac{\pi}{2}$. Clearly, the set of all such $(x,y,z)$ lies on the straight line segment joining $P$ and $Q$. In Figure 2, it is marked by the blue line $PQ$ where the plane $\Sigma$ cuts through the cube with side length $\frac{\pi}{2}$. Similarly, the two other blue line segments $PR$ and $QR$ correspond to right-angled triangles for which $y = \frac{\pi}{2}$ and $x = \frac{\pi}{2}$, respectively. In other words, all the points on the blue-lined triangle $PQR$, except the points $P,Q$ and $R$, correspond to different equivalence classes of right-angled triangles. 

\bigskip

\noindent
{\bf Acute-angled triangles:-} In this case, we have $0 < x,y,z < \frac{\pi}{2}$. The inside portion of the cube drawn in Figure 2 consists of all such points. We only therefore need to find which of these points lie on the plane $\Sigma$. Those are precisely the points lying strictly inside the triangular region $PQR$ except the line segments $PQ$, $QR$ and $RP$. 

\bigskip

\noindent
{\bf Obtuse-angled triangles:-} Since any triangle is either acute-angled or right-angled or obtuse-angled, we only need to find the portion on the plane $\Sigma$ that has not yet been discussed above. And this is precisely the region on $\Sigma$ that lies outside the cube. We can also verify directly that the points in the region correspond to obtuse-angled triangles. For, if a triangle is obtuse-angled, then $x > \frac{\pi}{2}$ or $y > \frac{\pi}{2}$ or $z > \frac{\pi}{2}$. Clearly, this represents the portion of $\Sigma$ that lies outside the cube.

\bigskip

\noindent
{\bf Isosceles triangles:-} In Figure $1$ of Section \ref{ds-section}, we have seen that the points on the line segments $BC$, $CD$ and $DE$ except the points $B, C$ and $E$ represent different classes of isosceles triangles depending on whether the two equal sides are bigger or smaller than the third side. 

\smallskip

Let us address this issue for our representation in Figure $2$ of Section \ref{amader}. We first note that an arbitrary point on the line segment $AR$, barring $A$ and $R$, is given by $t(0,0,\pi) + (1 - t)(\frac{\pi}{2},\frac{\pi}{2},0) = ((1 - t)\frac{\pi}{2},(1 - t)\frac{\pi}{2},t\pi)$ for some $t \in (0,1)$. Since two co-ordinates are equal for any such point, we conclude that these represent equivalence classes of isosceles triangles. In particular, for $t = \frac{1}{3}$, we obtain the equivalence class of equilateral triangles. Since the line segment $AR^{\prime}$ lies outside the cube, the points on this except $A$ and $R^{\prime}$ represent equivalence classes of obtuse-angled isosceles triangles. Following the same reasoning, we can say that the points on $RR^{\prime}$ except $R$ and $R^{\prime}$ represent equivalence classes of acute-angled isosceles triangles. We further note that the points on the line segment joining $A$ and $(\frac{\pi}{3},\frac{\pi}{3},\frac{\pi}{3})$, except the end points, represent the equivalence classes of those isosceles triangles for which the lengths of the equal sides are smaller than that of the third side. 

\smallskip

We can also run the same argument without taking into consideration the ordering of the side lengths of triangles. In that case the points on $CQ$ and $BP$, barring the points $B$, $C$, $P$ and $Q$, represent equivalence classes of isosceles triangles and the portions on these line segments lying outside (respectively, inside) represent obtuse-angled (respectively, acute-angled) isosceles triangles. Lastly, not surprisingly, we conclude by noting that the three straight lines $AR$, $BP$ and $CQ$ are the medians of the triangle $ABC$ and they are concurrent at the centroid $(\frac{\pi}{3},\frac{\pi}{3},\frac{\pi}{3})$ which represents the equivalence class of equilateral triangles.

\section{Concluding remarks: an attempt towards a  quantitative analysis}\label{amader-2}

Our method enables one to represent different class of similar triangles in $\mathbb{R}^{3}$ via the map described above. Under this correspondence, we see that the equivalence class of right-angled triangles are given by points on straight lines, those of acute-angled triangles are given by points inside a triangular region and those of obtuse-angled triangles are given by the points inside the union of three triangles. If one wishes to find some other type of triangles, for example, right-angled isosceles triangles, then we narrow down our focus on the boundary of the triangle $\Delta PQR$ as it corresponds to the equivalence classes of right-angled triangles. Now, the angles of a right-angled isosceles triangle are $\frac{\pi}{4},\frac{\pi}{4}$ and $\frac{\pi}{2}$. Therefore, the points $P^{\prime},Q^{\prime}$ and $R^{\prime}$ represent the equivalence class of right-angled isosceles triangles. Unlike in the treatment by Dhar and Sinha here we obtain three points instead of one because we do not take the orderings of the angles into consideration.

\smallskip

We now wish to understand the proportion of different types of equivalence classes of triangles in $\mathcal{T}$. For the sake of brevity, we denote by $\mathcal{T}_{o}$, $\mathcal{T}_{r}$ and $\mathcal{T}_{a}$ the equivalence classes of obtuse-angled, right-angled and acute-angled triangles, respectively. Since $\mathcal{T}$, $\mathcal{T}_{o}$, $\mathcal{T}_{r}$ and $\mathcal{T}_{a}$ are all infinite sets, the natural question is how does one define the notion of measurement of the sizes of these sets. The map $\Psi$ from $\mathcal{T}$ to the portion of $\Sigma$ lying in the first octant offers a way to accomplish this.

\smallskip

We observe that $\Psi(\mathcal{T})$ is the open triangular region $ABC$, which is a bounded subset of $\mathbb{R}^{3}$, as depicted in Figure 2. Moreover, $\Psi(\mathcal{T}_{o})$ is the union of the three open triangular regions $APQ$, $CPR$ and $BQR$. Since the area of these three triangular regions is $\frac{3}{4}$-th of the area of the triangular region $ABC$, we may reasonably declare that the proportion of $\mathcal{T}_{o}$ in $\mathcal{T}$ is $\frac{3}{4}$. Similarly, $\Psi(\mathcal{T}_{a})$ is the open triangular region $PQR$ which has area $\frac{1}{4}$-th of that of the triangular region $ABC$. In the same manner, we can declare that the proportion of $\mathcal{T}_{a}$ in $\mathcal{T}$ is $\frac{1}{4}$. Lastly, the equivalence classes of right-angled triangles are represented by the points on the perimeter of the triangle $PQR$ which has area $0$ in $\mathbb{R}^{3}$. Thus the proportion of $\mathcal{T}_{r}$ in $\mathcal{T}$ is obtained as $0$. 

\smallskip

Thus we see that the map $\Psi$ enables us to invoke the notion of measurement of subsets of $\mathcal{T}$ via the area of its image under $\Psi$. With this notion of measurement, we see that the proportion of the equivalence classes of obtuse-angled triangles is more than that of acute-angled and right-angled triangles. Also, it shows that the equivalence classes of right-angled triangles are very thin, even though there are infinitely many inequivalent right-angled triangles. This is comparable with a result in number theory that states that the proportion of prime numbers in the set of positive integers is $0$ even though there are infinitely many prime numbers. This is a consequence of the prime number theorem. The notion of proportion in number theory is captured by {\it density} (The enthusiastic reader is encouraged to look into the materials of \cite{tene} to learn more about densities). In our case, it is captured by the area of a region in the Euclidean space. 

\smallskip

One might wonder at this point as to what conclusions we can draw if we follow this approach of computing area in the exposition of Dhar and Sinha. As noted in Section \ref{ds-section}, in Figure 1, the points inside the region enclosed by the line segments $BE$, $ED$ and the arc $BD$, together with all the points on the line segment $ED$ except $E$ and $D$, represent all the equivalence classes of obtuse-angled triangle. The area of this region is $\frac{\pi - 2}{8}$. Again, the points on and inside the region bounded by the line segments $DC$, $CB$ and the arc $BD$ (except the points on the arc $BD$) represent all the equivalence classes of acute-angled triangles. The area of this region is $\frac{4 - \pi}{8}$. Therefore, here also we obtain that the proportion of obtuse-angled triangles is bigger than the proportion of acute-angled triangles.

\bigskip

\begin{ack}
It is our great pleasure to thank Dr. Shailesh Shirali for meticulously going through the manuscript and inspiring us to come up with some consequences of our construction made in Section \ref{amader}. This results in the quantitative analysis done in Section \ref{amader-2}. 
\end{ack}

\end{document}